\documentclass[12pt,a4paper]{amsart}

\usepackage {amsfonts}
\usepackage[english]{babel}
\def\b{\beta}
\def\g{\gamma}
\def\G{\Gamma}
\def\d{\delta}
\def\a{\alpha}
\def\b{\beta}
\def\p{\varphi}
\def\e{\varepsilon}
\def\l{\lambda}
\def\L{\Lambda}
\def\t{\theta}
\def\R{{\mathbb R}}
\def\C{{\mathbb C}}
\def\N{{\mathbb N}}
\def\Z{{\mathbb Z}}

\def\o{\omega}
\def\p{\varphi}
\def\bs{~\hfill\rule{7pt}{7pt}}
\def\la{\langle}
\def\ra{\rangle}

\DeclareMathOperator{\supp}{supp}
\DeclareMathOperator{\ord}{ord}

\newtheorem{Th}{Theorem}
\newtheorem{Cor}{Corollary}

\newtheorem{Pro}{Proposition}

\theoremstyle{definition}

\begin{document}

\title{Tempered distributions with discrete support and spectrum}

\author{S.Yu. Favorov}

\address{Sergii Favorov,
\newline\hphantom{iii}  Karazin's Kharkiv National University
\newline\hphantom{iii} Svobody sq., 4,
\newline\hphantom{iii} 61022, Kharkiv, Ukraine}
\email{sfavorov@gmail.com}

\maketitle {\small
\begin{quote}
\noindent{\bf Abstract.}
We investigate properties of tempered distributions with discrete or countable supports such that their Fourier transforms are distributions with discrete or countable supports as well. We find sufficient conditions for support of the distribution to be a finite union of translations of a full-rank lattice.
\medskip


\end{quote}
}
          \section{introduction}

Denote by $S(\R^d)$ Schwartz space of test functions $\p\in C^\infty(\R^d)$ with finite norms
 \begin{equation}\label{n}
  N_{n,m}(\p)=\sup_{\R^d}\{\max\{1,|x|^n\}\max_{\|k\|\le m} |D^k\p(x)|\},\quad n,m=0,1,2,\dots,
 \end{equation}
where
$$
k=(k_1,\dots,k_d)\in(\N\cup\{0\})^d,\ \|k\|=k_1+\dots+k_d,\  D^k=\partial^{k_1}_{x_1}\dots\partial^{k_d}_{x_d}.
 $$
 These norms generate the topology on $S(\R^d)$.  Elements of the space $S^*(\R^d)$ of continuous linear functionals on $S(\R^d)$ are called tempered distributions. For each tempered distribution $f$ there are $c<\infty$ and $n,\,m\in\N\cup\{0\}$ such that for all $\p\in S(\R^d)$
\begin{equation}\label{d}
                           |f(\p)|\le cN_{n,m}(\p).
\end{equation}
Moreover, this estimate is sufficient for distribution $f$ to belong to $S^*(\R^d)$
(see \cite{V}, Ch.3).

The Fourier transform of a tempered distribution $f$ is defined by the equality
$$
\hat f(\p)=f(\hat\p)\quad\mbox{for all}\quad\p\in S(\R^d),
$$
where
$$
   \hat\p(y)=\int_{\R^d}\p(x)\exp\{-2\pi i\la x,y\ra\}dx
 $$
is the Fourier transform of the function $\p$. Note that the Fourier transform of every tempered distribution is also a tempered distribution.

An element $f\in\R^d$ is called {\it a quasicrystal Fourier} if $f$ and $\hat f$ are discrete measures on $\R^d$. In this case the support of $\hat f$ is called {\it spectrum} of $f$. These notions were inspired by experimental discovery made in the middle of 80's of non-periodic atomic structures with diffraction patterns consisting of spots. There are a lot of papers devoted to investigation of properties of quasicrystals Fourier (see, for example, collections of works \cite{D}, \cite{Q}, papers \cite{F1}-\cite{M}, and so on).

 A set $A\subset\R^d$ is {\it discrete} if it has no finite limit points, and $A$ is {\it uniformly discrete} if it has a strictly positive  separating constant
 $$
 \eta(A):=\inf\{|x-x'|:\,x,\,x'\in A,\,x\neq x'\}.
  $$
  Complex Radon measure or distribution is discrete (uniformly discrete) if its support is discrete (uniformly discrete). we will call the support of Fourier transform of a tempered distribution $f$  spectrum of $f$ as well. Following \cite{L2}, we will say that a discrete set $A\subset\R^p$ is  {\it a set of finite type}, if the set
  $$
  A-A=\{x-x':\,x,x'\in A\}
  $$
  is discrete. A set $L\subset\R^d$ is called {\it a full-rank lattice} if  $L=T\Z^d$ for some nondegenerate linear operator $T$ on $\R^d$, the lattice $L^*=(T^t)^{-1}\Z^d$ is called the {\it conjugate lattice} for $L$. A set $A$ is a {\it pure crystal} with respect to a full-rank lattice $L$ if it is a finite union of  cosets of $L$.

\medskip
We begin with the following result of N.Lev and A.Olevskii \cite{LO1}, \cite{LO2} on quasicrystals.
\begin{Th}\label{TO}
Let $\mu$ be a uniformly discrete positive quasicrystal Fourier on $\R^d$ with uniformly discrete spectrum. Then the support of $\mu$ is a subset of a pure crystal with respect to a full-rank lattice $L$, and $\hat\mu$ is a subset of a pure crystal with respect to the conjugate lattice $L^*$. The same assertion is valid for an arbitrary uniformly discrete quasicrystal Fourier with the spectrum of a finite type. In the dimension $d=1$ the assertion is valid for every uniformly discrete quasicrystal Fourier with uniformly discrete spectrum. On the other hand, there are quasicrystals Fourier with discrete support and spectrum such that the above assertions do not valid.
\end{Th}
Note that for dimension $d>1$ there are non-positive quasicrystals Fourier with uniformly discrete support and spectrum such that their support is not a pure crystal (\cite{F2}).

There is another type of results.
\begin{Th}[F.\cite{F3}]\label{TF}
Let $\mu$ be a complex measure on $\R^d$ with discrete support $\L$ of a finite type such that $\inf_{x\in\L}|\mu(x)|>0$, let the Fourier transform be a measure $\hat\mu=\sum_{y\in\G} b(y)\d_y$ with the countable $\G\subset\R^d$ such that
$$
\sum_{|y|<r}|b(y)|=O(r^T),\quad r\to\infty,
$$
with some $T<\infty$. Then $\L$ is a pure crystal.
\end{Th}

Let us now go over from measures to distributions. An analog of Theorem \ref{TO} for tempered distributions was proved by V.Palamodov \cite{P}.
\begin{Th}
Let $f\in S^*(\R^d)$ be such that its support $\L$ and spectrum $\G$ are discrete sets of a finite type and, moreover, one of the differences $\L-\L$ and $\G-\G$ is uniformly discrete. Then $\L$ is a pure crystal with respect to a lattice $L$ and $\G$ is a pure crystal with respect to the conjugate lattice $L^*$.
\end{Th}
  In the present paper we obtain two analogs of Theorem \ref{TF} for tempered distributions.

     \section{the main results}

By \cite{Ru}, every distribution $f$ with discrete support $\L$ has the form
\begin{equation}\label{r}
f=\sum_{\l\in\L} P_\l(D)\d_\l, \quad P_\l(x)=\sum_{\|k\|\le K_\l}p_{\l,k}x^k, \quad x\in\R^d,\ p_{\l,k}\in\C,\ K_\l<\infty.
\end{equation}
 Here $\d_y$ means, as usual, the unit mass at the point $y\in\R^d$ and $x^k=x_1^{k_1}\cdots x_d^{k_d}$.

 Moreover, $\ord f=\sup_\l \deg P_\l<\infty$ (see Proposition \ref{P1} below).
Therefore we will consider distributions
\begin{equation}\label{r1}
f=\sum_{\l\in\L}\sum_{\|k\|\le m}p_{\l,k}D^k\d_\l,\quad k\in(\N\cup\{0\})^d.
\end{equation}
 If the Fourier transform $\hat f$ has a discrete support $\G$, we also have
\begin{equation}\label{r2}
\hat f=\sum_{\g\in\G}\sum_{\|j\|\le m'}q_{\g,j}D^j\d_\g,\quad j\in(\N\cup\{0\})^d.
\end{equation}
We will suppose that $m=\ord f$ and $m'=\ord\hat f$. Also, we will consider the case of distributions $f$ and $\hat f$ of forms (\ref{r1}) and (\ref{r2}) with arbitrary countable $\L$ and $\G$. If this is the case, we will also say that $\L$ is support and $\G$ is spectrum of $f$.

Denote by $B(x,r)$ the ball in $\R^d$ of radius $r$ with the center in $x$, $B(r)=B(0,r)$, by $\#A$  denote a number of elements of a discrete set $A$,  and put $n_A(r)=\#(A\cap B(r))$. We say that discrete $A$  is of {\it bounded density} if
$$
\sup_{x\in\R^d}\#(A\cap B(x,1))<\infty.
 $$
 Clearly, every uniformly discrete set $A$ is of  bounded density, and every set $A$  of bounded density satisfies the condition
 $$
 n_A(r)=O(r^d),\quad r\to\infty.
 $$
 A set $A\subset\R^d$ is {\it relatively dense} if there is $R<\infty$ such that every ball of radius $R$ intersects with $A$.

Also, for any $f$ of the form (\ref{r1}) set
$$
\kappa_f(\l)=\sup_{\|k\|\le m}|p_{\l,k}|,\qquad \rho_f(r)=\sum_{|\l|<r}\kappa_{f}(\l).
$$
Using  properties of almost periodic measures and sets, we prove the following theorems in Section 4.

\begin{Th}\label{T1}
Let $f_1,f_2$ be tempered distributions on $\R^d$ with discrete supports $\L_1,\,\L_2$, respectively, such that $\L_1-\L_2$ is a discrete set and
  \begin{equation}\label{k1}
  \inf_{\l\in\L_j}\kappa_{f_j}(\l)>0,\quad j=1,2.
  \end{equation}
If $\hat f_1,\,\hat f_2$ are both measures with countable supports such that
$$
\rho_{\hat f_1}(r)+\rho_{\hat f_2}(r)=O(r^T),\quad r\to\infty,
$$
 with some $T<\infty$,  then $\L_1,\L_2$ are pure crystals with respect to a unique full-rank lattice.
 \end{Th}

\begin{Cor}\label{C1}
 Let $f$ be a tempered distribution on $\R^d$ with discrete support $\L$ of a finite type such that  $\inf_{\l\in\L}\kappa_{f}(\l)>0$.
If $\hat f$ is a measure with countable support such that
$$
\rho_{\hat f}(r)=O(r^T),\quad r\to\infty
$$
with some $T<\infty$, then $\L$ is a pure crystal.
 \end{Cor}

\begin{Th}\label{T2}
 Let $f_1,f_2$ be tempered distributions on $\R^d$ with discrete relatively dense supports $\L_1,\,\L_2$ and discrete spectrums $\G_1,\,\G_2$, let $\L_1-\L_2$ be a discrete set,  and let
 $$
 n_{\G_1}(r)+n_{\G_2}(r)=O(r^T),\quad r\to\infty
 $$
 with some $T<\infty$. If conditions (\ref{k1}) and
  \begin{equation}\label{k2}
   \sup_{\l\in\L_j}\kappa_{f_j}(\l)<\infty,\quad j=1,2,
 \end{equation}
 satisfy,  then $\L_1,\L_2$ are pure crystals with respect to a unique full-rank lattice.
 \end{Th}

\begin{Cor}\label{C2}
 Let $f$ be tempered distribution on $\R^d$ with discrete support $\L$ of a finite type and discrete spectrum $\G$ such that
 $$
 n_\G(r)=O(r^T),\quad r\to\infty
 $$
  with some $T<\infty$. If
 $$
   0<\inf_{\l\in\L}\kappa_{f}(\l)\le\sup_{\l\in\L}\kappa_{f}(\l)<\infty,
  $$
  then $\L$ is a pure crystal.
 \end{Cor}
\section{preliminary properties of distributions with discrete supports}

\begin{Pro}\label{P1}
 Suppose $f\in S^*(\R^d)$ has form (\ref{r}) with  discrete $\L$ and satisfies (\ref{d}) for some $n,m$. Then
 $$
 \sup_{\l\in\L}\deg P_\l\le m,
 $$
 and there exists $C<\infty$ such that
 \begin{equation}\label{p}
  |p_{\l,k}|\le C\max\{1,|\l|^n\} \mbox{ for all $\l\in\L$ and $k$ such that }\|k\|=m.
 \end{equation}
 If $\L$ is uniformly discrete, then for all $k$, $\|k\|\le m$,
 \begin{equation}\label{e}
 |p_{\l,k}|\le C\max\{1,|\l|^n\}.
 \end{equation}
 \end{Pro}
 The second part of the Proposition was earlier proved by V.Palamodov \cite{P}.
\medskip

{\bf Proof of Proposition \ref{P1}}. Set $\l\in\L$ and  $\e\in(0,1)$ such that
$$
\inf\{|\l-\l'|:\,\l'\in\L,\,\l'\neq\l\}>\e.
 $$
 Let $\p$ be a function on $\R$ such that
 $$
 \p(|x|)\in C^\infty(\R^d),\quad\p(|x|)=0\mbox{ for }|x|>1/2,\quad \p(|x|)=1\mbox{ for }|x|<1/3.
  $$
  Put
  $$
  \p_{\l,k,\e}(x)=\frac{(x-\l)^k}{k!}\p(|x-\l|/\e)\in S(\R^d),
  $$
   where, as usual, $k!=k_1!\cdots k_d!$. It is easily shown that
   $$
   f(\p_{\l,k,\e})=(-1)^{\|k\|}p_{\l,k}.
   $$
  On the other hand, we get for some $c(\a,\b)<\infty$, where $\a,\b\in(\N\cup\{0\})^d$,
$$
 |f(\p_{\l,k,\e})|\le\sup_{|x-\l|<\e}\max\{1,|x|^n\}\sum_{\|\a+\b\|\le m} c(\a,\b)\left|D^\a\p\left(\frac{|x-\l|}{\e}\right)D^\b\left(\frac{(x-\l)^k}{k!}\right)\right|.
$$
Note that
$$
|D^\a\p(|x-\l|/\e)|\le \e^{-\|\a\|}c(\a)
$$
 and
$$
D^\b(x-\l)^k=\left\{\begin{array}{l} 0\mbox{ if }k_j<\b_j\mbox{ for at least one }j,\\c(k,\b)(x-\l)^{k-\b}\mbox{ if }k_j\ge\b_j\mbox{ for all }j .
  \end{array}\right.
$$
Since $|x-\l|<1$, we get
$$
\max\{1,|x|^n\}\le 2^n\max\{1,|\l|^n\}.
$$
Taking into account that
$$
\supp\p(|x-\l|/\e)\subset B(\l,\e),
$$
 we get
\begin{equation}\label{p1}
 |p_{\l,k}|\le\sum_{\|\a+\b\|\le m,\,\b_j\le k_j\,\forall j} c(k,\a,\b)\max\{1,|\l|^n\}\e^{\|k\|-\|\a+\b\|}.
\end{equation}
Note that coefficients $c(k,\a,\b)$ do not depend on $\l$ and $\e$.

 If $\|k\|>m$,  we take $\e\to0$ and obtain $p_{\l,k}=0$, hence $\deg P_\l\le m$ for all $\l$.

If $\|k\|=m$, we obtain $|p_{\l,k}|\le C\max\{1,|\l|^n\}$.

If $\L$ is uniformly discrete, we pick
$$
\e=\e_0<\min\{1,\,\eta(\L)/2\},
$$
 and (\ref{p1}) implies the estimate
 $$
 |p_{\l,k}|\le C\e_0^{-m}\max\{1,|\l|^n\}
 $$
 for all $k$. \bs
\medskip

 {\bf Remark}. There are tempered distributions with discrete support such that conditions (\ref{e}) do not valid for $\|k\|<m$.

    Set
 $$
 \l_j=j,\ \l'_j=j+2^{-2j},\ \L=\{\l_j,\,\l'_j\}_{j\in\N},\ f=\sum_{j\in\N}2^j(\d_{\l'_j}-\d_{\l_j}).
 $$
  For any $\psi\in S(\R)$ we get
 $$
 f(\psi)=\sum_{j\in\N}2^j(\psi(j+2^{-2j})-\psi(j))=\sum_{j\in\N}2^{-j}\psi'(j+\t_j2^{-2j})
 $$
 with $\t_j\in (0,1)$. Therefore,
 $$
 |f(\psi)|\le\sup_{x\in\R}|\psi'(x)|,
  $$
  and $f\in S^*(\R)$. Since
  $$
  p_{\l'_j,0}=-p_{\l_j,0}=2^j,
  $$
 we see that estimate (\ref{e}) does not valid for $k=0$.
\medskip

\begin{Pro}\label{P2}
  Suppose $f\in S^*(\R^d)$ satisfies (\ref{d}) with some $n,\,m$, $f$ has form (\ref{r1}) with countable $\L$, and $\hat f$ has form (\ref{r2}) with discrete $\G$. Then
  $$
  n\ge d,\quad \ord\hat f\le n-d,\quad |q_{\g,j}|\le C'\max\{1,|\g|^m\}\mbox{ for }\|j\|=n-d.
  $$
  For the case of uniformly discrete $\G$ we get
  $$
  |q_{\g,j}|\le C\max\{1,|\g|^m\}\quad\forall j,\  \|j\|\le n-d.
  $$
  \end{Pro}

\begin{Cor}\label{C3}
 Suppose $f\in S^*(\R^d)$ satisfies (\ref{d}) with $n=d$ and some $m$, has form (\ref{r1}) with countable $\L$, and spectrum $\G$ of $f$ is discrete. Then $\hat f$ is a measure, and
 $$
 \hat f=\sum_{\g\in\G}q_\g\d_\g,\quad |q_\g|\le C'\max\{1,|\g|\}^m.
 $$
\end{Cor}

{\bf Proof of Proposition \ref{P2}}.  Set $\g\in\G$ and pick $\e\in(0,1)$ such that
$$
\inf\{|\g-\g'|:\,\g'\in\G,\,\g'\neq\g\}>\e.
 $$
 Let $\p$ be the same as in the proof of Proposition \ref{P1}. Put
 $$
 \p_{\g,l,\e}(y)=\frac{(y-\g)^l}{l!}\p(|y-\g|/\e)\in S(\R^d).
 $$
  We have
 \begin{equation}\label{q}
 (-1)^{\|l\|}q_{\g,l}=\sum_{\|j\|\le m'}q_{\g,j}D^j\d_\g(\p_{\g,l,\e}(y))=(\hat f,\p_{\g,l,\e})=(f,\hat\p_{\g,l,\e}).
 \end{equation}
Note that
 $$
 \hat\p_{\g,l,\e}(x)=e^{-2\pi i\la x,\g\ra}(l!)^{-1}(-2\pi i)^{-\|l\|}D^l(\widehat{\p(\cdot/\e)})=c(l)e^{-2\pi i\la x,\g\ra}\e^{d+\|l\|}(D^l\hat\p)(\e x).
 $$
 Therefore,
 $$
 D^k(\hat\p_{\g,l,\e}(x))=\e^{d+\|l\|}\sum_{\a+\b=k}c(\a,\b)D^\a\left[e^{-2\pi i\la x,\g\ra}\right]D^\b[(D^l\hat\p)(\e x)]
 $$
 $$
 =\sum_{\a+\b=k}c(\a,\b)(-2\pi i)^{\|\a\|}\g^\a e^{-2\pi i\la x,\g\ra}\e^{d+\|l\|+\|\b\|}(D^{\b+l}\hat\p)(\e x).
 $$
 Since $\hat\p\in S(\R^d)$, we get for every $k$, $\|k\|\le m$, and every $M<\infty$
  \begin{equation}\label{D}
 |D^k(\hat\p_{\g,l,\e}(x))|\le C(M)\e^{d+\|l\|}\max\{1,|\g|^m\}(\max\{1,|\e x|^M\})^{-1}.
 \end{equation}
 Pick $M=n$. By (\ref{d}), we have
 $$
    |(f,\hat\p_{\g,l,\e})|\le c(f)\sup_{\R^d}[\max\{1,|x|^n\}\max_{\|k\|\le m} |D^k(\hat\p_{\g,l,\e}(x))|].
 $$
 Hence,
 $$
  |q_{\g,l}|=|(f,\hat\p_{\g,l,\e})|\le C'(\hat\p)\e^{d+\|l\|}\max\{1,|\g|^m\}\sup_{x\in\R^d}[\max\{1,|x|^n\}
  (\max\{1,|\e x|^n\})^{-1}].
  $$
 Since $\sup$ in the right-hand side of the inequality equals $\e^{-n}$, we obtain
 $$
 |q_{\g,l}|\le C'\max\{1,|\g|^m\}\e^{\|l\|+d-n}.
 $$
 If $\|l\|>n-d$, we take $\e\to0$ and get $q_{\g,l}=0$, hence $\ord\hat f\le n-d$.

 For $\|l\|=n-d$ we get $|q_{\g,l}|\le C'\max\{1,|\g|^m\}$.

 If $\G$ is uniformly discrete,  we take $\e=\e_0<\eta(\G)/2$ for all $\g\in\G$
  and obtain the bound
  $$
  |q_{\g,l}|\le \e_0^{d-n} C'\max\{1,|\g|^m\}\quad\forall l,\ \|l\|\le n-d.\qquad\qquad\qquad  \bs
  $$

\medskip
{\bf Remark}. Since
$$
N_{n,m}(\hat\p)\le C(n,m) N_{m,n}(\p)
$$
 for every $\p\in S(\R^d)$, we see that estimate (\ref{d}) for $f\in S^*(\R^d)$ implies the estimate $|\hat f(\p)|\le c'N_{m,n}(\p)$. Therefore the direct application of Proposition \ref{P1} gives only the bound $\ord\hat f\le n$.

\medskip
 M.Kolountzakis, J.Lagarias proved in \cite{KL} that Fourier transform  of every measure $\mu$ on the line $\R$ with  support of bounded density, bounded masses $\mu(x)$, and discrete spectrum is also a measure $\hat\mu=\sum_{\g\in\G}q_\g\d_\g$ with uniformly bounded $q_\g$. The following proposition generalizes this result for distributions from $S^*(\R^d)$.

 \begin{Pro}\label{P3}
Suppose $f\in S^*(\R^d)$  has form (\ref{r1}) with some $m$ and countable $\L$, and discrete spectrum $\G$. If
 $$
\rho_f(r)=O(r^{d+H}),\quad r\to\infty,\quad H\ge0,
  $$
  then $\ord\hat f\le H$. If, in addition, $H$ is integer, then for $\|j\|=H$ we get $|q_{\g,j}|\le C'\max\{1,|\g|^m\}$. For the case of uniformly discrete $\G$ we get
  $$
  |q_{\g,j}|\le C\max\{1,|\g|^m\}\quad\forall j,\quad \|j\|\le H.
  $$
\end{Pro}
\begin{Cor}\label{C4}
 If $f\in S^*(\R^d)$  has form (\ref{r1}) with some $m$ and countable $\L$, discrete spectrum $\G$,  and
 $$
 \rho_f(r)=O(r^d)\quad r\to\infty,
 $$
 then $\hat f$ is a measure, and
 $$
 \hat f=\sum_{\g\in\G}q_\g\d_\g,\ |q_\g|\le C'\max\{1,|\g|^m\}.
 $$
\end{Cor}

\begin{Cor}\label{C4a}
 If $f\in S^*(\R^d)$  has form (\ref{r1}) with some $m$, bounded $\kappa_f(\l)$, discrete support $\L$ of bounded density, and discrete spectrum $\G$, then $\hat f$ is a measure, and
 $$
 \hat f=\sum_{\g\in\G}q_\g\d_\g,\ |q_\g|\le C'\max\{1,|\g|^m\}.
 $$
\end{Cor}

 {\bf Proof of Proposition \ref{P3}}. Let $\p_{\g,l,\e}$ be the same as in the proof of Proposition \ref{P2}. By (\ref{q}) and (\ref{r1}),
  $$
   (-1)^{\|l\|}q_{\g,l}=(f,\hat\p_{\g,l,\e})=\sum_{\l\in\L}\sum_{\|k\|\le m}p_{\l,k}(-1)^{\|k\|}D^k(\hat\p_{\g,l,\e}(\l)).
  $$
  Using (\ref{D}), we get
  $$
  |q_{\g,l}|\le C'(m,M)\e^{d+\|l\|}\max\{1,|\g|^m\}\sum_{\l\in\L}\kappa_f(\l)(\max\{1,|\e\l|^M\})^{-1}.
  $$
We have
 $$
 \sum_{\l\in\L}\kappa_f(\l)(\max\{1,|\e\l|^M\})^{-1}=\rho_f(1/\e)+\e^{-M}\int_{1/\e}^\infty t^{-M}d\rho_f(t)
  $$
 Pick $M>d+H$.   Integrating by parts and using the estimate for $\rho_f(r)$, we see that the right-hand side is equal to
 $$
O(\e^{-d-H})+\e^{-M}M\int_{1/\e}^\infty \rho_f(t)t^{-M-1}dt=O(\e^{-d-H})\quad\hbox{ as }\e\to\infty.
 $$
 Finally,
 $$
 |q_{\g,l}|\le C'\max\{1,|\g|^m\}\e^{\|l\|-H}.
 $$
  If $\|l\|>H$, we take $\e\to0$ and get $q_{\g,l}=0$, hence, $m'=\ord\hat f\le H$.

 If $H$ is integer, we get $|q_{\g,l}|\le C'\max\{1,|\g|^m\}$ for $\|l\|=H$.

 If $\G$ is uniformly discrete,  we take $\e=\e_0<\eta(\G)/2$ for all $\g\in\G$ and obtain the bound
  $$
 |q_{\g,l}|\le \e_0^{-H} C'\max\{1,|\g|^m\}\quad \forall l, \|l\|\le m'.\qquad\qquad\qquad  \bs
 $$

\section{almost periodic distributions and proofs of Theorems \ref{T1} and \ref{T2}}

Recall that a continuous function $g$ on $\R^d$ is  almost periodic if for any  $\e>0$ the set of $\e$-almost periods of $g$
  $$
  \{\tau\in\R^d:\,\sup_{x\in\R^d}|g(x+\tau)-g(x)|<\e\}
  $$
is a relatively dense set in $\R^d$.

Almost periodic functions are uniformly bounded  on $\R^d$. The class of almost periodic functions is closed with respecting to taken absolute values, linear combinations, maximum, minimum of a finite family of almost periodic functions. A limit of a uniformly in $\R^d$ convergent sequence  of almost periodic functions is also almost periodic. A typical example of an almost periodic function is an absolutely convergence exponential sum $\sum c_n\exp\{2\pi i\la x,\o_n\ra\}$ with $\o_n\in\R^d$ (see, for example, \cite{C}).

A measure $\mu$ on $\R^d$ is called almost periodic if the function
 $$
 (\psi\star\mu)(t)=\int_{\R^d}\psi(x-t)d\mu(x)
 $$
  is almost periodic in $t\in\R^d$ for each continuous function $\psi$ on $\R^d$ with compact support. A distribution $f\in S^*(\R)$ is almost periodic if the function $(\psi\star f)(t)=f(\psi(\cdot-t))$ is almost periodic in $t\in\R^d$ for each $\psi\in S(\R^d)$ (see \cite{A},\cite{M},\cite{R}). Clearly, every almost periodic distribution has a relatively dense support.

\begin{Th}\label{T6}
Let $f\in S^*(\R^d)$ and $\hat f$ be a measure $\sum_{\g\in\G}q_\g\d_\g$ with countable $\G$ such that
\begin{equation}\label{ro}
\rho_{\hat f}(r)=O(r^T),\quad r\to\infty
\end{equation}
for some $T<\infty$. Then $f$ is an almost periodic distribution.
\end{Th}

{\bf Proof of Theorem \ref{T6}}. Let $\psi\in S(\R^d)$. The converse Fourier transform of the function $\psi(x-t)$ is $\hat\psi(-y)e^{2\pi i\la y,t\ra}$.
Therefore,
 \begin{equation}\label{s}
 f(\psi(\cdot-t))=\hat f(\hat\psi(-y)e^{2\pi i\la y,t\ra})=\sum_{\g\in\G}q_\g\hat\psi(-\g)e^{2\pi i\la\g,t\ra}.
 \end{equation}
 Taking into account that $\hat\psi(-y)\in S(\R^d)$, we see that the sum in (\ref{s}) is majorized by the sum
$$
\sum_{\g\in\G}C(\psi,T)|q_\g|\max\{1,|\g|\}^{-T-1}=C(\psi,T)\left[\rho_{\hat f}(1)+\int_1^\infty r^{-T-1}d\rho_{\hat f}(r)\right].
$$
 Integrating by parts, we obtain that the integral converges,  then the  sum in (\ref{s}) converges absolutely, and $\psi\star f$ is almost periodic. \bs

Combining this theorem with the results of the previous section, we obtain
\begin{Th}\label{T7}
Suppose $f\in S^*(\R^d)$ has form (\ref{r1}) with countable $\L$, has discrete spectrum $\G$ such that
$n_\G(r)=O(r^T)$ as $r\to\infty$  for some $T<\infty$, and satisfies one of the following conditions

i) inequality (\ref{d}) holds with $n=d$ and some $m$;

ii)  $\rho_f(r)=O(r^d)$ as $r\to\infty$;

iii)  $\kappa_f(\l)$ is bounded and support $\L$ is of bounded density.
\smallskip

\noindent Then $f$ is an almost periodic distribution.
\end{Th}

{\bf Proof of Theorem \ref{T7}}. Using Corollaries \ref{C3}, or \ref{C4}, or \ref{C4a}, we obtain that $\hat f$ is a measure,
$$
\hat f=\sum_{\g\in\G}q_\g\d_\g,\qquad |q_\g|\le C'\max\{1,|\g|\}^m.
$$
Therefore, integrating by parts, we get
 $$
\rho_{\hat f}(r)\le C'\left[\rho_{\hat f}(1)+\int_1^r t^m dn_\G(t)\right]=O(r^{m+T}),\quad r\to\infty.
 $$
Then we apply Theorem \ref{T6} completes the proof. \bs
\medskip

Now we can proof the main results of our paper.
\smallskip

{\bf Proof of Theorem \ref{T1}}. By Theorem \ref{T6}, $f_j,\,j=1,2,$ are almost periodic distributions. In particular, $\L_1$ is relatively dense. If there are $\l_n, \l'_n\in\L_2$ such that $\l_n\neq\l'_n$ and $\l_n-\l'_n\to0$ as $n\to\infty$, then there are $R<\infty$ and $x_n\in\L_1$ such that $\l_n, \l'_n\in B(x_n,R)$. Therefore there exists infinitely many points of the set $\L_1-\L_2$ in the ball $B(R)$, that is impossible. Hence $\eta(\L_2)>0$ and, similarly, $\eta(\L_1)>0$. By Proposition \ref{P1}, $\ord f_j<\infty,\,j=1,2$.

We apply Proposition \ref{P3} to the measure $\tilde\mu(y)=\hat f_1(-y)$ with $m=0$ and $H=T-d$ (we may suppose that $T\ge d$). Since the Fourier transform of $\tilde\mu(y)$ is just the tempered distribution $f_1$ with uniformly discrete support, we get representation (\ref{r1}) for $f_1$ with
$$
\sup_{\l\in\L_1}\sup_{\|k\|\le\ord f_1}|p_{\l,k}|<\infty.
 $$
Put $s=\inf_{\l\in\L_1}\kappa_{f_1}(\l)$, and pick a number $\e\in(0,\eta(\L_1)/5)$
such that
$$
\sum_{\|k\|\le\ord f_1}|p_{\l,k}|<s/(2\e)\quad \forall \l\in\L_1.
$$
Let $\p_\e(x)$ be nonnegative $C^\infty$-function such that
 $$
 \supp\p_\e\subset B(\e+\e^2),\quad \p_\e(x)=1 \mbox{ for }|x|<\e.
 $$
  Put
  $$
  \p_{k,\e}(x)=\p_\e(x)x^k/k!,\quad k\in(\N\cup\{0\})^d.
   $$
  Since $f_1$ is almost periodic, we see that the function $g_{k,\e}=f_1\star\p_{k,\e}$ is almost periodic.

 Fix $\l\in\L_1$. By (\ref{k1}), there exists $k'$ such that $|p_{\l,k'}|\ge s$. Put
 $$
 A=\{k\in(\N\cup\{0\})^d:\,k_j\le k_j'\ \forall j=1,\dots,d\}.
  $$
 For $x\in B(\l,\e)$ we have
$$
g_{k',\e}(x)=\left(\sum_{\|k\|\le\ord f_1}p_{\l,k}D^k\d_\l\right)\frac{(\l-x)^{k'}}{{k'}!}=\sum_{k\in A}p_{\l,k}(-1)^{\|k\|}(\l-x)^{k'-k}/(k'-k)!
$$
 Therefore,
$$
|g_{k',\e}(x)|\ge |p_{\l,k'}|-\e\sum_{k\in A,k\neq k'}|p_{\l,k}|>s/2.
$$
Now set
$$
h_\e(x)=\min\{1,2s^{-1}\sup_{\|k\|\le\ord f_1}|g_{k,\e}(x)|\}.
$$
 Clearly, $h_\e(x)$ is an almost periodic function and
$$
h_\e(x)=1\mbox{ for } x\in\bigcup_{\l\in\L_1}B(\l,\e),\qquad \supp h_\e\subset\bigcup_{\l\in\L_1}B(\l,\e+\e^2),\qquad 0\le h_\e(x)\le1.
$$
Let $\psi$ be an arbitrary continuous function on $\R^d$ with support in the ball $B(\eta(\L_1)/5)$. It is readily seen that
the function $\psi\star h_\e$ is almost periodic as well. Since $\e<\eta(\L_1)/5$, we see that for each fixed $t\in\R^d$ the support of the function $\psi(\cdot-t)$ intersects with at most one ball $B(\l,\e+\e^2)$. Therefore we get
 \begin{eqnarray}\label{h1}
  \left|\frac{\psi\star h_\e(t)}{\o_d\e^d}-\frac{1}{\o_d\e^d}\int_{B(\l,\e)}h_\e(x)\psi(x-t)dx\right|
  \le \left|\frac{1}{\o_d\e^d}\int_{B(\l,\e+\e^2)\setminus B(\l,\e)}h_\e(x)\psi(x-t)dx\right|\nonumber \\                                                       \le\frac{\o_d(\e+\e^2)^d-\o_d\e^d}{\o_d\e^d}\sup_{\R^d}|\psi(x)|,\qquad\qquad\qquad\qquad\qquad\qquad,
 \end{eqnarray}
where $\o_d$ is a volume of the unit ball in $\R^d$, and
 \begin{equation}\label{h2}
  \left|\frac{1}{\o_d\e^d}\int_{B(\l,\e)}h_\e(x)\psi(x-t)dx-\psi(\l-t)\right|\le\sup_{x\in B(\l,\e)}|\psi(x-t)-\psi(\l-t)|.
   \end{equation}
It follows from (\ref{h1}) and (\ref{h2}) that the almost periodic functions $(\o_d\e^d)^{-1}(\psi\star h_\e)$ uniformly converge as $\e\to0$ to the function $\psi\star\d_{\L_1}$, where $\d_{\L_1}=\sum_{\l\in\L_1}\d_{\l}$. Hence the function $\psi\star\d_{\L_1}$ is almost periodic. It is readily seen that the same is true for every continuous function $\psi$ on $\R^d$ with compact support. The similar construction works for $\d_{\L_2}=\sum_{\l\in\L_2}\d_{\l}$.

Thus  $\d_{\L_1}, \d_{\L_2}$ are almost periodic measures with discrete set of differences $\L_1-\L_2$. Applying Theorem 6 from \cite{F2}, we obtain that $\L_1,\L_2$ are pure crystals with respect to a unique full-rank lattice. \bs

\medskip
{\bf Proof of Theorem \ref{T2}}. Since the sets $\L_j,\,j=1,2$, are relatively dense, we get, arguing as before,  the bounds $\eta(\L_j)>0$ for $j=1,2$. Hence, both $\L_j$ are of bounded density. By Corollary \ref{C4a}, $\hat f_j$ are measures such that
$$
\kappa_{\hat f_j}(\l)=O(|\l|^m),\quad |\l|\to\infty,
$$
for some $m<\infty$. Using Theorem \ref{T7} iii), we obtain that both $f_j$ are almost periodic distributions, hence we can repeat the proof of Theorem \ref{T1}.\bs
\medskip

{\bf Remark}. It is clear that the bound $\eta(\L)>0$ follows immediately from the conditions of Corollaries \ref{C1} or \ref{C2}.

\end{document}